\documentclass[11pt]{article}
\usepackage{amssymb,amsmath}
\usepackage{float}
\usepackage{color,fullpage}

\usepackage{algorithm,algorithmic}
\usepackage{amsthm}
\newcommand{\RR}{\mathbb{R}}

\newcommand{\cX}{{\mathcal{X}}}

\DeclareMathOperator*{\Min}{minimize}

\newtheorem{theorem}{Theorem}
\newtheorem{lemma}{Lemma}

\begin{document}

\title{Linear Convergence of the Proximal Incremental Aggregated Gradient Method under Quadratic Growth Condition}

\author{Hui Zhang\thanks{
College of Science, National University of Defense Technology,
Changsha, Hunan, 410073, P.R.China. Corresponding author. Email: \texttt{h.zhang1984@163.com}}
}



\date{\today}

\maketitle

\begin{abstract}
 Under the strongly convex assumption, several recent works studied the global linear convergence rate of the proximal incremental aggregated gradient (PIAG) method for minimizing the sum of a large number of smooth component functions and a non-smooth convex function. In this paper, under \textsl{the quadratic growth condition}--a strictly weaker condition than the strongly convex assumption, we derive a new global linear convergence rate result, which implies that the PIAG method attains global linear convergence rates in both the function value and iterate point errors. The main idea behind is to construct a certain Lyapunov function.
\end{abstract}

\newline



\section{Introduction}
A fundamental generic optimization model arises in many problems in machine learning, signal processing, image science, communication systems, and distributed optimization to mention just a few. It consists in minimizing the sum of a differentiable function $F(x)$ and a possibly non-smooth regularization function $h(x)$:
\begin{equation}\label{main}
\Min_{x\in \RR^d} \Phi (x)\triangleq F(x)+h(x).
\end{equation}
The well-known method to solve this model is the forward-backward splitting (FBS) scheme, which is often called the proximal gradient method in the literature. This method consists of the composition of a gradient (forward) step of $F(x)$ with a proximal (backward) step on $h(x)$. It can be expressed as follows:
\begin{align}
& y_k=x_k-\alpha \cdot \nabla F(x_k), \\
& x_{k+1}=\arg\min_x\{h(x)+\frac{1}{2\alpha}\|x-y_k\|^2\},
\end{align}
where $\alpha>0$ is some step size. The remarkable merit of this method lies in exploiting the smooth plus smooth structure of model \eqref{main}. However, if the differentiable function $F(x)$ is a sum of $N$ component functions $f_n(x)$, i.e.,
\begin{equation}\label{Fx}
F(x)=\sum_{n=1}^N f_n(x)
\end{equation}
with large $N$ (this indeed appears in many applications), then evaluating the full gradient of $F(x)$, i.e., $\nabla F(x)=\sum_{n=1}^N \nabla f_n(x)$ is costly and even prohibitive. Hence, a natural idea to overcome this difficulty is to modify the standard FBS method by utilizing the additional structure \eqref{Fx}. Following this line of thought, the proximal incremental aggregated gradient (PIAG) method was proposed and studied in several recent papers \cite{Global2016Vanli,Analysis2016Aytekin,A2016Vanli}. At each iteration $k\geq 0$, it first constructs a vector that aggregates the gradients of all components functions, possibly evaluated at $k-\tau_k^n$ iteration,
$$g_k=\sum_{n=1}^N\nabla f_n(x_{k-\tau_k^n})$$
where $\tau_k^n$ are some nonnegative integers.  This vector is used to exploit the additive structure \eqref{Fx} of the $N$ component functions, but also to approximate the full gradient of $F(x)$ at the current iteration. In fact, if $\tau_k^n\equiv 0$, then $g_k=\nabla F(x_k)$. After obtaining $g_k$, the PIAG method then performs a proximal step on the sum of the non-smooth term $h(x)$ and $\langle g_k, x-x_k\rangle$ as follows:
\begin{equation}\label{xk}
x_{k+1}=\arg\min_x\{h(x)+\langle g_k, x-x_k\rangle + \frac{1}{2\alpha}\|x-x_k\|^2\}.
\end{equation}
By introducing an auxiliary vector, we can rewrite the PIAG method into the following scheme:
\begin{align}
& y_k=x_k-\alpha \cdot g_k, \\
& x_{k+1}=\arg\min_x\{h(x)+\frac{1}{2\alpha}\|x-y_k\|^2\}.
\end{align}
We can see that the PIAG method differs from its mother scheme--the FBS method mainly at the gradient step, and reduces to it as  $\tau_k^n\equiv 0$.

In this paper, we will focus on the global linear convergence of the PIAG method. Beforehand, we review several very recent works around this topic.

\subsection{Related work}
Work \cite{Global2016Vanli} is the first study that establishes a global linear convergence rate for the PIAG method in function value error, i.e., $\Phi(x_k)-\Phi(x^*)$, where $x^*$ denotes the minimizer point of $\Phi(x)$.  Work \cite{Analysis2016Aytekin} used a different analysis and showed a global linear convergence rate in iterate point error, i.e., $\|x_k-x^*\|$. The authors of \cite{Global2016Vanli}  combined the results presented in \cite{Global2016Vanli} and \cite{Analysis2016Aytekin} and provided a stronger linear convergence rate for the PIAG method in the recent paper \cite{A2016Vanli}. However, all these mentioned works are built on the strongly convex assumption, which is actually not satisfied by many application problems and hence motives lots of research to find weaker alternatives. Influential weaker conditions include the error bound property, the restricted strongly convex property, the quadratic growth condition, and the Polyak-{\L}ojasiewicz inequality; the interested reader could refer to \cite{Zhang2015The,zhang2016new,I2015Linear,Karimil2016linear,Drusvyatskiy2016Error,Bolte2015From,Zhou2015A}. Works \cite{Zhang2015The,Drusvyatskiy2016Error}  studied the linear convergence of the FBS method under these weaker conditions. But to our knowledge, there is no work of studying the global linear convergence of the PIAG method under these weaker conditions.

\subsection{Main contribution}
Our main contribution of this paper is a global linear convergence result of the PIAG method under the quadratic growth condition, which is strictly weaker than the strongly convex assumption.

Although we employ an important lemma presented in \cite{Analysis2016Aytekin}, which was also used in \cite{A2016Vanli}, our proof strategy  is essentially different from that in \cite{Analysis2016Aytekin} and \cite{A2016Vanli}. The main idea behind the proof strategy is to construct a Lyapunov function that simultaneously includes the function value and iterate point errors. Hence, our linear convergence result simultaneously characterizes the global linear convergence rate in both the function value and iterate point errors; see Theorem 1 in Section 3. This perspective is new even for the FBS method.

Besides, the global linear convergence rate derived in this study has a linear dependence on the condition number of the problem and a quadratic dependence on the delay parameter. This is consistent with the results in \cite{Analysis2016Aytekin} and \cite{Global2016Vanli}.  Moreover, if the delay parameter is less than $47$, then we recover the result from \cite{A2016Vanli}, which has the best (linear) dependence on the condition number and the delay parameter. Although the constraint of the delay parameter being less than $47$  is not required in \cite{A2016Vanli}, there used the strongly convex assumption and the analysis is essentially more complicated than ours.

\section{Assumptions}First, we list the following standard assumptions that are used in both \cite{Analysis2016Aytekin} and \cite{A2016Vanli}.
\begin{itemize}
  \item[A1.] Each component function $f_n(x)$ is convex with $L_n$-continuous gradient, that is,
  $$\|\nabla f_n(x)-\nabla f_n(y)\|\leq L_n\|x-y\|, \forall x, y\in \RR^d.$$
  This assumption implies that the sum function $F(x)$ is convex with $L$-continuous gradient, where $L=\sum_{n=1}^N L_n$.
  \item[A2.] The regularization function $h(x):\RR^d\rightarrow (-\infty, \infty]$ is proper, closed, convex and subdifferentiable everywhere in its effective domain, i.e., $\partial h(x)\neq \emptyset$ for all $x\in \{y\in\RR^d: h(x)<\infty\}$.
  \item[A3.] The time-varying delays $\tau_k^n$ are bounded, i.e., there is a nonnegative integer $\tau$ such that
  $$\tau_k^n\in \{0,1, \cdots, \tau\},$$
  hold for all $k\geq 1$ and $n\in\{1,2,\cdots, N\}$. Such $\tau$ is called the delay parameter.
\end{itemize}
Besides the standard assumptions A1-A3, we replace the strongly convex condition used previously by the following quadratic growth condition:
\begin{itemize}
  \item[A4.] The objective function $\Phi(x)$ satisfies the quadratic growth condition, meaning there is a real number $\beta>0$ such that
      $$\Phi(x)-\Phi^*\geq \frac{\beta}{2}d^2(x,\cX),$$
      where $\cX$ is the set of  minimizers of $\Phi(x)$, which is assumed to be nonempty, $\Phi^*$ is the minimal value of $\Phi(x)$, and $d(x,\cX)$ is the distance function from points to a fixed set, defined by
$$d(x,\cX)=\inf_{y\in \cX}\|x-y\|.$$
\end{itemize}

On one hand, the strong convexity with parameter $\mu>0$ implies that
$$\Phi(x)-\Phi^*\geq \frac{\mu}{2}\|x-x^*\|^2,$$
where $x^*$ is the unique minimizer of $\Phi(x)$, and hence implies the quadratic growth condition. On the other hand, we can easily construct functions to show that the quadratic growth condition does not imply any strong convexity. For example, the composition $g(Ax)$, where $g(\cdot)$ is a strongly convex and $A$ is rank deficient, satisfies the quadratic growth condition but fails to be strongly convex. Therefore, the quadratic growth condition is \textsl{strictly weaker} than the strongly convex condition. The former has been adopted recently as an efficient alternative of the latter to derive linear convergence rate results for many fundamental algorithms, such as the proximal gradient method \cite{Drusvyatskiy2016Error,Zhang2015The}, the conditional gradient method \cite{Beck2015Linearly}, and the cyclic block coordinate gradient descent method \cite{zhang2016new}. Due to the close connection of the quadratic growth condition with previously mentioned weaker notions \cite{Zhang2015The,zhang2016new,I2015Linear,Karimil2016linear,Drusvyatskiy2016Error,Bolte2015From,Zhou2015A},  the methods in \cite{Zhou2015A,Li2016Calculus} can be used to verify whether a given function satisfies the quadratic growth condition.

\section{Main results}
Throughout this section, we remind the reader that, we consider the optimization model \eqref{main} with $F(x)$ given by \eqref{Fx} and the sequence $\{x_k\}$ generated by the PIAG method.

First, we introduce the following lemma, which was presented in \cite{Analysis2016Aytekin}.
\begin{lemma}
Assume that the nonnegative sequences $\{V_k\}$ and $\{w_k\}$ satisfy the following inequality:
$$V_{k+1}\leq a V_k -b w_k+ c\sum_{j=k-k_0}^k w_j,$$
for some real numbers $a\in(0,1)$ and $b,c\geq 0$, and some nonnegative integer $k_0$. Assume also that $w_k=0$ for $k<0$, and the following holds:
$$\frac{c}{1-a}\frac{1-a^{k_0+1}}{a^{k_0}}\leq b.$$
Then, $V_k\leq a^k V_0$ for all $k\geq 0$.
\end{lemma}

Before presenting the main result of this paper, we state another lemma, which can be viewed as a generalization of the standard result for the proximal gradient method; see lemma 2.3 in \cite{beck2009fast}.
\begin{lemma}
Suppose that the standard assumptions A1-A3 hold. Let
$$\Delta_k=\frac{L(\tau+1)}{2}\sum_{j=k-\tau}^k\|x_{j+1}-x_j\|^2.$$
Then, the following holds:
\begin{equation}
\Phi(x_{k+1})\leq \Phi(x)+\frac{1}{2\alpha}\|x-x_k\|^2-\frac{1}{2\alpha}\|x-x_{k+1}\|^2-\frac{1}{2\alpha}\|x_{k+1}-x_k\|^2+\Delta_k.
\end{equation}
\end{lemma}

Now, we present the main theorem of this paper.
\begin{theorem}
Suppose that the assumptions A1-A4 hold, and the step-size satisfies:
$$\alpha\leq \frac{\left(1+\frac{\beta}{L}\frac{1}{\tau+1}\right)^{\frac{1}{(\tau+1)}}-1}{\beta}.$$
Define a Lyapunov function $$\Psi(x)\triangleq\Phi(x)-\Phi^*+\frac{1}{2\alpha}d^2(x,\cX).$$
Then, the PIAG method is linearly convergent in the sense that
\begin{equation}\label{result1}
\Psi(x_k)\leq \left(1-\frac{\alpha\beta}{1+\alpha\beta}\right)^k \Psi(x_0),
\end{equation}
for all $k\geq0$. In particular, the PIAG method attains a global linear convergence in function value error:
\begin{equation}\label{result2}
\Phi(x_k)-\Phi^* \leq \left(1-\frac{\alpha\beta}{1+\alpha\beta}\right)^k \Psi(x_0),
\end{equation}
and a global linear convergence in iterate point error:
\begin{equation}\label{result3}
d^2(x_k,\cX)\leq \Psi(x_0) \frac{2\alpha }{1+\alpha\beta}\left(1-\frac{\alpha\beta}{1+\alpha\beta}\right)^k ,
\end{equation}
for all $k\geq0$. Furthermore, if
$$\alpha=\frac{\left(1+\frac{\beta}{L}\frac{1}{\tau+1}\right)^{\frac{1}{(\tau+1)}}-1}{\beta},$$
then
\begin{equation}\label{result4}
\Psi(x_k)\leq \left(1-\frac{1}{(\tau+1)(\tau+2)\eta}\right)^k \Psi(x_0),
\end{equation}
for all $k\geq0$, where $\eta=L/\beta$ stands for the number condition of optimization problem \eqref{main}.
\end{theorem}

Some comments are in order:

\begin{itemize}
  \item First, since the strongly convex assumption implies the quadratic growth condition with $\beta=\mu$, we have that
$$1-\frac{\alpha\beta}{1+\alpha\beta}=\frac{1}{1+\alpha\beta}=\frac{1}{\mu\alpha+1},$$
and hence the global linear convergence in iterate point error \eqref{result3} attains the following form:
\begin{equation}\label{result3}
d^2(x_k,\cX)\leq \Psi(x_0) \frac{2\alpha }{1+\alpha\mu}\left(\frac{1}{\mu\alpha+1}\right)^k.
\end{equation}
This recovers the main result from \cite{Analysis2016Aytekin}.
  \item Second, if the delay parameter $\tau \leq 47$, then following from \eqref{result4} we have the following linear convergence in function value error:
  \begin{equation}
\Phi(x_k)\leq \left(1-\frac{1}{49 \eta(\tau+1)}\right)^k \Psi(x_0).
\end{equation}
This recovers the main result from \cite{A2016Vanli}. Note that the constraint of $\tau \leq 47$ is not required in \cite{A2016Vanli}. But there used the strongly convex conditions as an assumption.

\end{itemize}

\section{Proofs}
\subsection{Proof of Lemma 2}
We divide the proof into two parts. The first part can be found from the proof of Theorem 1 in \cite{Analysis2016Aytekin}; we include it here for completion. The second part is a standard argument, which is different from the optimality condition based method adopted in the proof of Theorem in \cite{Analysis2016Aytekin}.

\textbf{Part 1.} Since each component function $f_n(x)$ is convex with $L_n$-continuous gradient, we have the following upper bound estimations:
\begin{align}
f_n(x_{k+1})\leq &f_n(x_{k-\tau_k^n})+\langle \nabla f_n(x_{k-\tau_k^n}), x_{k+1}-x_{k-\tau_k^n}\rangle +\frac{L_n}{2}\|x_{k+1}-x_{k-\tau_k^n}\|^2 \nonumber\\
\leq  &f_n(x)+\langle \nabla f_n(x_{k-\tau_k^n}), x_{k+1}-x\rangle +\frac{L_n}{2}\|x_{k+1}-x_{k-\tau_k^n}\|^2 \label{Lp}.
\end{align}
Summing \eqref{Lp} over all components functions and using the expression of $g_k$, we obtain
\begin{equation}
F(x_{k+1})\leq F(x)+\langle g_k, x_{k+1}-x\rangle +\sum_{n+1}^N\frac{L_n}{2}\|x_{k+1}-x_{k-\tau_k^n}\|^2.
\end{equation}
The last term of the inequality above can be upper-bounded using Jensen's inequality as follows:
 \begin{align}
 \sum_{n=1}^N\frac{L_n}{2}\|x_{k+1}-x_{k-\tau_k^n}\|^2=&\sum_{n=1}^N\frac{L_n}{2}\|\sum_{j=k-\tau_k^n}^k (x_{j+1}-x_j)\|^2 \nonumber \\
 \leq &\frac{L(\tau+1)}{2}\sum_{j=k-\tau}^k\|x_{j+1}-x_j\|^2=\Delta_k.
 \end{align}
 Therefore,
 \begin{equation}\label{LP1}
F(x_{k+1})\leq F(x)+\langle g_k, x_{k+1}-x\rangle + \Delta_k.
\end{equation}

 \textbf{Part 2.} By the definition of $x_{k+1}$, $x_{k+1}$ is the minimizer of the $\frac{1}{\alpha}$-strongly convex function
 $$x\mapsto h(x)+\langle g_k, x-x_k\rangle + \frac{1}{2\alpha}\|x-x_k\|^2,$$
 hence for all $x\in \RR^d$, we have
 \begin{equation}\label{SC}
h(x_{k+1})\leq h(x)+\langle g_k, x-x_{k+1} \rangle + \frac{1}{2\alpha}\|x-x_k\|^2-\frac{1}{2\alpha}\|x-x_{k+1}\|^2-\frac{1}{2\alpha}\|x_{k+1}-x_k\|^2.
\end{equation}
Adding \eqref{LP1} and \eqref{SC} yields the desired inequality.

\subsection{Proof of Theorem 1}
Below, we use $x^\prime$ to stand for the projection of $x$ onto the set $\cX$.
Let us write successively the inequality in Lemma 2 at $x=x_k$ and then at $x=x_k^\prime$. Note that $\Phi(x_k^\prime)=\Phi^*$. We obtain
\begin{equation}\label{ineq1}
\Phi(x_{k+1})\leq \Phi(x_k)-\frac{1}{\alpha}\|x_k-x_{k+1}\|^2+\Delta_k
\end{equation}
and
\begin{equation}\label{ineq2}
\Phi(x_{k+1})\leq \Phi^* +\frac{1}{2\alpha}\|x_k^\prime-x_k\|^2-\frac{1}{2\alpha}\|x_k^\prime-x_{k+1}\|^2-\frac{1}{2\alpha}\|x_{k+1}-x_k\|^2+\Delta_k.
\end{equation}
By the definition of projection, we have
\begin{equation}\label{ineq3}
 \|x_k^\prime-x_{k+1}\|^2\geq \|x_{k+1}^\prime-x_{k+1}\|^2=d^2(x_{k+1}, \cX).
\end{equation}
We split the term $\frac{1}{2\alpha}\|x_k^\prime-x_k\|^2$ into two terms as follows:
$$\frac{1}{2\alpha}\|x_k^\prime-x_k\|^2=\mu \|x_k^\prime-x_k\|^2+\nu\|x_k^\prime-x_k\|^2,$$
where $\mu+\nu=\frac{1}{2\alpha}$ and $\mu,\nu>0$. By the quadratic growth condition, we relax the second term to introduce the function value error $\Phi(x_k)-\Phi^*$. Thus,
\begin{equation}\label{ineq4}
\nu\|x_k^\prime-x_k\|^2=\nu d^2(x_{k}, \cX)\leq \frac{2\nu}{\beta}(\Phi(x_k)-\Phi^*).
\end{equation}
Now, using \eqref{ineq2} together with \eqref{ineq3} and \eqref{ineq4}, we obtain the following relation
\begin{equation}\label{ineq5}
\Phi(x_{k+1})-\Phi^* +\frac{1}{2\alpha}d^2(x_{k+1}, \cX) \leq \frac{2\nu}{\beta}(\Phi(x_k)-\Phi^*) +\mu  d^2(x_{k}, \cX)-\frac{1}{2\alpha}\|x_{k+1}-x_k\|^2+\Delta_k.
\end{equation}
As we will show later, this relation is sufficient for us to derive the desired result. But we now proceed in a slightly complicated way to include possible benefit brought by the relation \eqref{ineq1}.

 Multiplying the relation \eqref{ineq1} by a parameter $\lambda\geq 0$ and adding the resulting inequality and \eqref{ineq5}, we obtain
\begin{align}\label{ineq6}
&\Phi(x_{k+1})-\Phi^* +\frac{1}{2\alpha(1+\lambda)}d^2(x_{k+1}, \cX) \nonumber  \\ \leq& \frac{\lambda+\frac{2\nu}{\beta}}{1+\lambda}\left(\Phi(x_k)-\Phi^* +\frac{\mu}{\lambda+\frac{2\nu}{\beta}}  d^2(x_{k}, \cX)\right)-\frac{\lambda+0.5}{ \alpha(1+\lambda)}\|x_{k+1}-x_k\|^2+\Delta_k.
\end{align}
Denote $a(\lambda,\nu)=\frac{\lambda+\frac{2\nu}{\beta}}{1+\lambda}$. In order to apply Lemma 1 to derive the desired linear convergence rate, we require $a(\lambda,\nu)\in (0,1)$, which can be guaranteed by letting $0<\nu<\frac{\beta}{2}$. Besides, we require $$\frac{\mu}{\lambda+\frac{2\nu}{\beta}} \leq \frac{1}{2\alpha(1+\lambda)}$$
to relax the first term on the right-hand side of \eqref{ineq6}. This can be guaranteed by letting
$$\nu\geq \frac{1}{\frac{2}{\beta}+2(1+\lambda)\alpha}\triangleq \nu_0.$$
Now, for any fixed $\lambda\geq 0$, we use this lower bound of $\nu$ to yield the smallest rate:
$$a(\lambda, \nu_0)=1-\frac{\alpha\beta}{1+\alpha\beta(1+\lambda)}.$$
Thus, the smallest rate for varying parameter $\lambda\geq 0$ is $1-\frac{\alpha\beta}{1+\alpha\beta}$. In other words, the relation \eqref{ineq1} can not help to improve the linear convergence rate. Now, we set $\lambda=0, \nu=\nu_0$ in \eqref{ineq6}. Note that for such setting,
$$\frac{\mu}{\lambda+\frac{2\nu}{\beta}}=\frac{1}{2\alpha},~~
\frac{\lambda+\frac{2\nu}{\beta}}{1+\lambda}=1-\frac{\alpha\beta}{1+\alpha\beta}\triangleq a.$$
Thus, we use \eqref{ineq6} to obtain
\begin{align}\label{ineq7}
&\Phi(x_{k+1})-\Phi^* +\frac{1}{2\alpha}d^2(x_{k+1}, \cX)\nonumber \\ \leq & a\left(\Phi(x_k)-\Phi^* +\frac{1}{2\alpha}d^2(x_k, \cX)\right)-\frac{1}{ 2\alpha }\|x_{k+1}-x_k\|^2+\Delta_k,
\end{align}
i.e.,
\begin{equation}\label{ineq8}
\Psi(x_{k+1})  \leq   a\Psi(x_k)-\frac{1}{ 2\alpha }\|x_{k+1}-x_k\|^2+\frac{L(\tau+1)}{2}\sum_{j=k-\tau}^k\|x_{j+1}-x_j\|^2,
\end{equation}
where we use the expressions of $\Psi(x)$ and $\Delta_k$.
Note that $\|x_{j+1}-x_j\|^2=0$ for all $j<0$. We are ready to apply Lemma 1 with $V_k=\Psi_k$, $w_k=\|x_{k=1}-x_k\|^2$, $a=1-\frac{\alpha\beta}{1+\alpha\beta}$, $b=\frac{1}{ 2\alpha }$, $c=\frac{L(\tau+1)}{2}$, and $k_0=\tau$. To ensure Lemma 1 hold, we need
$$\frac{c}{\alpha\beta(1+\alpha\beta)^{-1}}\frac{1-(1+\alpha\beta)^{-\tau-1}}{(1+\alpha\beta)^{-\tau}}\leq \frac{1}{2\alpha}$$
to hold. Simplifying and rearranging terms we obtain the desired upper bound for the step-size $\alpha$. The linear convergence result \eqref{result1} follows from Lemma 1. In particular, \eqref{result2} is a direct consequence of \eqref{result1}, and \eqref{result3} follows from the quadratic growth condition and \eqref{result1}.

Finally, it remains to show \eqref{result4}. Taking the certain value $\alpha=\frac{\left(1+\frac{\beta}{L}\frac{1}{\tau+1}\right)^{\frac{1}{(\tau+1)}}-1}{\beta}$ in \eqref{result1} and using $\eta=\frac{L}{\beta}$, we derive that
\begin{align}
\Psi(x_k)& \leq\left(1+\frac{1}{\eta(\tau+1)}\right)^{\frac{-k}{\tau+1}} \Psi(x_0) \nonumber \\
&\leq \left(1-\frac{1}{\eta(\tau+2)}\right)^{\frac{k}{\tau+1}} \Psi(x_0) \nonumber \\
&\leq \left(1-\frac{1}{\eta(\tau+2)(\tau+1)}\right)^k \Psi(x_0),
\end{align}
where the second line follows as $\eta\geq 1$, and in the third line we use the Bernoulli inequality, i.e., $(1+x)^r\leq 1+rx$ for any $x\geq -1$ and $r\in [0,1]$.  This completes the proof.

\section*{Acknowledgements}
This work is supported by the National Science Foundation of China (No.11501569).  The main idea of this research was carried out
while the author was visiting BeiJing International Center for Mathematical Research by invitation of Professor Zaiwen, Wen.

\bibliographystyle{abbrv}


\end{document}